\documentclass[a4paper,12pt,reqno]{amsart}

\usepackage[english]{babel}
\usepackage[T1]{fontenc}
\usepackage[left=1in,right=1in,top=1.2in,bottom=1.2in]{geometry}
\usepackage{times}
\usepackage{microtype}

\usepackage{commath,amssymb,amscd,mathrsfs,mathtools,dsfont,bbm,multirow,array,caption,comment,pifont}
\usepackage[all]{xy}
\usepackage{enumerate}
\usepackage{longtable}

\usepackage{setspace}
\allowdisplaybreaks

\usepackage[dvipsnames,svgnames,table]{xcolor}
\definecolor{red}{RGB}{255,25,25}
\definecolor{blue}{RGB}{25,50,200}
\usepackage{tikz-cd}
\usepackage[pagebackref,linktocpage]{hyperref}

\hypersetup{
pdftitle={},				
pdfauthor={Hu--Truong--Xie},	
colorlinks=true,			
linkcolor=red,			
citecolor=MidnightBlue,	
filecolor=magenta,		
urlcolor=MidnightBlue	
}

\usepackage{cleveref}	

\newtheorem{theorem}{Theorem}[section]
\crefname{theorem}{Theorem}{Theorems}
\newtheorem{lemma}[theorem]{Lemma}
\crefname{lemma}{Lemma}{Lemmas}
\newtheorem{proposition}[theorem]{Proposition}
\crefname{proposition}{Proposition}{Propositions}

\crefname{prop}{Proposition}{Propositions}
\newtheorem{corollary}[theorem]{Corollary}
\crefname{corollary}{Corollary}{Corollaries}

\crefname{cor}{Corollary}{Corollaries}

\crefname{conjecture}{Conjecture}{Conjectures}

\crefname{conj}{Conjecture}{Conjectures}
\newtheorem*{conj*}{Conjecture}
\crefname{conj}{Conjecture}{Conjectures}

\crefname{conjA}{Conjecture}{Conjecture}

\crefname{conjB}{Conjecture}{Conjecture}

\crefname{conjC}{Conjecture}{Conjecture}

\crefname{conjDk}{Conjecture}{Conjecture}

\crefname{conjD}{Conjecture}{Conjecture}

\crefname{conjH}{Conjecture}{Conjecture}

\crefname{conjGr}{Conjecture}{Conjecture}

\theoremstyle{definition}

\crefname{definition}{Definition}{Definitions}

\crefname{defn}{Definition}{Definitions}

\crefname{example}{Example}{Examples}

\crefname{notation}{Notation}{Notation}
\newtheorem*{notation*}{Notation}
\crefname{notation}{Notation}{Notation}

\crefname{problem}{Problem}{Problems}

\crefname{question}{Question}{Questions}

\crefname{condition}{Condition}{Conditions}

\crefname{assumption}{Assumption}{Assumptions}

\crefname{propGr}{Property}{Property}

\theoremstyle{remark}

\crefname{rmk}{Remark}{Remarks}
\newtheorem*{rmk*}{Remark}
\crefname{rmk}{Remark}{Remarks}

\crefname{remark}{Remark}{Remarks}

\crefname{fact}{Fact}{Facts}

\crefname{claim}{Claim}{Claims}
\newtheorem*{claim*}{Claim}
\crefname{claim}{Claim}{Claims}

\crefname{step}{Step}{Steps}

\crefname{case}{Case}{Cases}

\numberwithin{equation}{section}

\newcommand{\doi}[1]{\href{https://doi.org/#1}{{\tt DOI:#1}}}

\newcommand{\bQ}{\mathbf{Q}}
\newcommand{\bR}{\mathbf{R}}

\newcommand{\sT}{\mathsf{T}}

\newcommand{\id}{\operatorname{id}}

\newcommand{\N}{\mathsf{N}}

\newcommand{\Tr}{\operatorname{Tr}}

\newcommand{\Z}{\mathsf{Z}}

\begin{document}

\title[Tate's question, Standard conjecture D, semisimplicity and Dynamical degrees]{Tate's question, Standard conjecture D, semisimplicity  and Dynamical degree comparison conjecture}

\author{Fei Hu}
\address{School of Mathematics, Nanjing University, Nanjing 210093, China}

\email{\href{mailto:fhu@nju.edu.cn}{\tt fhu@nju.edu.cn}}

\author{Tuyen Trung Truong}
\address{Department of Mathematics, University of Oslo, 0316 Oslo, Norway}
\email{\href{mailto:tuyentt@math.uio.no}{\tt tuyentt@math.uio.no}}

\author{Junyi Xie}
\address{Beijing International center for Mathematical research, Peking University, China}
\email{\href{mailto:xiejunyi@bicmr.pku.edu.cn}{\tt xiejunyi@bicmr.pku.edu.cn}}

\begin{abstract}

Let $X$ be a smooth projective variety of dimension $n$ over the algebraic closure of a finite field $\mathbb{F}_p$.

Assuming the standard conjecture $D$, we prove
\begin{itemize}
    \item a weaker form of the Dynamical Degree Comparison conjecture;
    \item equivalence of semisimplicity of Frobenius endomorphism and of any polarized endomorphism (a more general result, in terms of the biggest size of Jordan blocks, holds).
\end{itemize}
We illustrate these results through examples, including varieties dominated by rational maps from Abelian varieties and suitable products of $K3$ surfaces.

Using the same idea, we provide a new proof of the main result in a recent paper by the third author, including Tate's question/Serre's conjecture that for a polarized endomorphism  $f:X\rightarrow X$, all eigenvalues of the action of $f$ on $H^k(X)$ have the same absolute value.

\end{abstract}


\date{\today}


\maketitle


\section{Introduction}

Let $\mathbf{K}$ be the algebraic closure of a finite field $\mathbb{F}_q$. We fix in the whole paper an irreducible smooth projective variety $X$ of dimension $n$ over $\mathbf{K}$. 

On $X$, there are so-called Weil cohomology theories $H^*(X)$, with good properties, similar to the singular cohomology groups on a compact K\"ahler manifold. Each group $H^k(X)$ can be (non-canonically) identified with a finitely dimensional complex vector space over the field of complex numbers $\mathbf{C}$. An example is the etale groups $H_{etal}^k(X,\mathbf{Q}_l)$, where $l$ is a prime number different from the characteristic of $\mathbf{K}$. For more detail, the reader can see \cite{KM74}. In this paper, we fix once and for all such a Weil cohomology theory $H^k(X)$. 

Our main objectives in this paper are the Semi-simplicity conjecture and Standard conjecture D, which will be recalled now, and the Weaker Dynamical degree comparison conjecture (to be discussed later after sufficient preliminary details). We recall that for $k=0,\ldots ,\dim (X)$, the numerical group $N^k(X)$ is the quotion of the group of algebraic cycles of pure dimension $k$ on $X$ modulo the numerical equivalence. The reader can consult \cite{Kleiman68} for more detail.  

{\bf Standard Conjecture D.} (Grothendieck and Bombieri)  $X$ satisfies the Standard conjecture D on the cohomology group $H^{2k}(X)$ if for every algebraic cycle $V$ of pure codimension $k$ in $X$, $V$ is $0$ in $N^k(X)$ iff the fundamental class $[V]$ is $0$ in $H^{2k}(X)$. 

{\bf Semisimplicity conjecture.} \cite[\S3, Conjecture (d)]{Tate65}  Let $f:X\rightarrow X$ be a polarized endomorphism. Then, for each $k$, the pullback of a polarized endomorphism  $F^*:H^*(X)\rightarrow H^*(X)$ is semisimple, i.e. diagonalisable.

{\bf Remark:} In case $f$ is an endomorphism such that there is an ample divisor $L$ so that $f^*(L)\sim L$, then the map $f$ has a very simple structure, and $f^*$ is semisimple. 

Here is a historical context  and motivation for the above conjectures and the Weaker Dynamical degree comparison conjecture. They started with the well known: 

{\bf Weil's Riemann hypothesis:} Let $\lambda$ be an eigenvalue of $Fr_q^*:H^k(X)\rightarrow H^k(X)$. Then $|\lambda |=q^{k/2}$. 

This conjecture inspired a lot of research in Algebraic Geometry since the 1960s. In particular, Grothendieck and Bombieri introduced the Standard conjectures (see \cite{Kleiman68}), modeled on a suggestion from \cite{Serre60}, whose validity will solve Weil's Riemann hypothesis. Surprisingly, Deligne \cite{Deligne74} solved Weil's Riemann hypothesis without solving the Standard conjectures. Interestingly, the proof of Deligne also yielded the proofs of several of the Standard conjectures. Nowadays, except the Standard conjecture C (which is a consequence of Deligne's proof, see \cite{KM74}) the other conjectures are largely unsolved: Lefschetz type standard conjecture (Standard conjecture A, B), Standard conjecture D and the Hodge Standard conjecture (a counterpart of the well known Hodge-Riemann relations on compact K\"ahler manifolds). If the Hodge standard conjecture holds, then Standard conjectuers A, B and D are equivalent, for many other relations between the Standard conjectures, the reader can consult \cite{Kleiman68}. In contrast to the situation over fields of zero characteristic, over finite fields the current status seems to indicate that the Hodge Standard conjecture is much harder than Standard conjectures A, B, D. In positive characteristic, besides easy examples (like products of projective spaces), Standard conjecture D is proven for Abelian varieties \cite{Clozel99} and certain self-products of K3 surfaces \cite{IIK22}, and the Hodge standard conjecture is known for Abelian varieties of dimension $4$ \cite{Ancona21} and certain self-products of K3 surfaces \cite{IIK22}.  

\cite{Serre60} proved Weil's Riemann hypothesis for polarized endomorphisms on compact K\"ahler manifolds. This inspired the Standard conjectures. Together with the Semisimplicity conjecture (which generalises Weil's Riemann hypothesis), the following conjecture, both are consequences of the validity of all Standard conjectures, are explicitly asked in \cite[\S3, Conjecture (d)]{Tate65}. 

{\bf Tate's question/Serre's conjecture.} Let $f:X\rightarrow X$ be a polarized endomorphism, i.e. there is an ample divisor $L$ and a positive integer $a>1$ so that $f^*(L)\sim aL$. Then, for each $k$ and each $\lambda $ an eigenvalue of $f^*:H^k(X)\rightarrow H^k(X)$, we have $|\lambda |=a^{k/2}$. 

Besides Abelian varieties case which was settled by Weil, Semisimplicity conjecture is still largely unkown until now. Tate's question/Serre's conjecture is recently solved by the third author \cite{Xie24}, which is an inspiration for the current paper. 

Now we come to the setup for Weaker Dynamical degree comparison conjecture. We consider a general morphism $f:X\rightarrow X$. Like cohomology groups, we can also pullback numerical groups by an endomorphism. We denote by $\chi _k(f)$ the spectral radius of $f^*:H^{k}(X)\rightarrow H^k(X)$ (for $k=0,1,\ldots, 2n$), and   $\lambda _k(f)$ the spectral radius of $f^*:N^k(X)\rightarrow N^k(X)$ (for $k=0,1,\ldots ,n$). The above conjectures can be naturally extended to this and even more general settings. 

We recall first the so-called {\bf Entropy comparison conjecture} from \cite{ES13}, inspired by the results on topological entropy of holomorphic endomorphisms of compact K\"ahler manifolds \cite{Gromov03,Yomdin87}, which amounts to showing the following equality for endomorphisms: $\max _{k=0,1,\ldots ,2n}\chi _k(f)=\max _{k=0,1,\ldots ,n}\lambda _k(f)$. \cite{ES13} itself solves the Entropy comparison conjecture for $X=$ a surface and $f$ an automorphism, where it is also mentioned a comment from Deligne that Entropy comparison conjecture should follow from the Standard conjectures. Entropy comparison conjecture is fully solved by the second author \cite{Truong}, and later an alternative proof is given in \cite{Shuddhodan19}.

Inspired results from complex dynamics (in particularly \cite{Gromov03,Yomdin87,DS05a,DS08,DN11,DNT12}), \cite{Truong} extended both Serre's conjecture and Entropy comparison conjecture to dynamical correspondences. A dynamical correspondence on $X$ is an effective algebraic cycle $V$ of codimension $n$ in $X\times X$, so that the natural projections to the first and second factors of each irreducible component $W$ are surjective. A dynamical correspondence is finite (developed by Suslin and Voevodsky, see \cite{MVW06}) if these natural projections are finite morphisms on each irreducible component. For two dynamical correspondences $f$ and $g$ on $X$, we can compose them dynamically $f\Diamond g$, to obtain a new dynamical correspondence. We refer the reader to \cite{DS08, Truong20, HT24} for an overview on this dynamical correspondence, where also the relation to the usual composition of algebraic cycles in Algebraic Geometry is discussed. To illustrate, if $f$ and $g$ are dominant rational maps, then $f\Diamond g$ is the same as the usual composition of dominant rational maps. If $f$ is a dynamical correspondence and $t$ is a positive integer, we denote by $f^{\Diamond t}=f\Diamond \ldots \Diamond f$ (t times) the t-th iterate. As endomorphisms, dynamical correspondences also pullback on cohomology and numerical groups. 

For a finitely dimensional complex vector space $H$, we choose a norm $||.||_H$ (and for convenience will write $||.||$ if no confusion arises). Since any two norms on $H$ are equivalence, for what discussed next a specific choice of $||.||$ does not make any difference. 

Let $f:X\vdash X$ be a dynamical correspondence. \cite{Truong20} shows that the following limit  $\lim _{t\rightarrow\infty}||(f^{\Diamond t})^*|_{N^k(X)}||^{1/t}$ exists, and hence we can denote it $\lambda _k(f)$ (which is a generalisation of the spectral radius of $f^*|_{N^k(X)}$ in case $f$ is an endomorphism).  An alternative approach, for rational maps $f$, is given in \cite{Dang20}. These are the counterparts of those on compact K\"ahler manifolds, established earlier in  \cite{DS05a, DS08}.  A sufficient condition to guarantee that the sequence $\lambda _j(f)$, $j=0,1,\ldots ,n$ is log-concave is that there is an infinite sequence $n_l$ such that the graphs of the iterates $f^{\Diamond n_l}$ are irreducible, see \cite[Theorem 1.1]{Truong20}. This is in particular the case if $f$ is a dominant rational map on $X$ (in this case, see also \cite{Dang20}). 

In \cite{Truong}, a generalisation of $\chi _k(f)$ to dynamical correspondences is also given. For $k=0,\ldots ,n$, we define: 
$$\chi _k(f)=\limsup _{t\rightarrow \infty}||(f^{\Diamond t})|_{H^k(X)}||^{1/t}.$$
It is unknown whether $\chi _k(f)$ is always finite, or if the limsup in the definition can be replaced by lim. On the other hand, \cite{Truong} shows that if Standard conjecture D holds, the natural extension of Entropy comparison conjecture to all dynamical correspondences also holds, and as a consequence we then have $\chi _k(f)<\infty$ for all $k$. 

The following two conjectures are for a dynamical correspondence $f:X\vdash X$. {\bf Dynamical degree comparison conjecture} (see \cite{Truong}): for all $k=0,\ldots ,n$, we have $\chi _{2k}(f)=\lambda _k(f)$, and moreover, $\chi _{2k+1}(f)\leq \sqrt{\lambda _k(f)\lambda _{k+1}(f)}$. {\bf Norm comparison conjecture} (see \cite{HT24}): There is a constant $C>0$ independent of $f$ such that for all $k=0,1,\ldots ,n$, we have $||f^*|_{H^{2k}(X)}||\leq C ||f^*|_{N^k(X)}||$, and $||f^*|_{H^{2k+1}(X)}||^2\leq C ||f^*|_{N^k(X)}||\times ||f^*|_{N^{k+1}(X)}||$. 

Dynamical degree comparison conjecture includes as a special case Serre's conjecture. For Abelian varieties and $f$ an endomorphism, it is solved by the first author \cite{Hu19,Hu24}. Still for Abelian varieties, but $f$ a general dynamical correspondence, Dynamical degree comparison conjecture is proved by the first two authors in \cite{HT24}. Norm comparison conjecture is yet more general than Dynamical degree comparison conjecture. 

Again, for  projective manifolds $X$ over $\mathbf{C}$, Dynamical degree comparison conjecture and Norm comparison conjecture hold. This follows from the fact that on $H^{2k}(X,\mathbf{R})$, we have a positivity property. More precisely, any smooth closed real $2k$-form is bounded by a multiple of $\omega ^k$, where $\omega$ has the same cohomological class as that of a very ample divisor. 

\cite{HT24} proposes an approach to resolve Norm comparison conjecture via the  use of a special cohomological correspondence $G_r$, defined as follows. For each rational number $r>0$, $G_r:H^*(X)\rightarrow H^*(X)$ is the  multiplication by $r^k$ on the vector space $H^k(X)$, for each $k=0,1,\ldots ,2n$. Over a finite field, as a consequence of Weil's Riemann hypothesis, Standard Conjecture C holds, and hence each $G_r$ is represented by an algebraic cycles. \cite{HT24} proposed a so-called Conjecture $G_r$, which stipulate some expected properties of this $G_r$ cohomological operator, which can be viewed as a quantified version of Standard Conjecture C. In Section 2 we will discuss more about this conjecture, including some new heuristic arguments on why (at least a version of it, enough for applications to the above conjectures) should hold. In the same paper, it is shown that if both Standard Conjecture D and Conjecture $G_r$ hold, the Norm comparison conjecture holds (and hence Dynamical degree comparison conjecture and Semisimplicity conjecture hold). In \cite{HT23}, it is shown that for polarized endomorphisms, Conjecture $G_r$ is a consequence of Standard conjectures.  

In the paper \cite{Xie24}, besides solving both the Tate's conjecture/Serre's question and more generally the Dynamical degree comparison for the case of endomorphisms, \cite{Xie24} also proved a weaker result for finite correspondences. The following is a generalisation of the arguments in \cite{Xie24}  to all dynamical correspondences.

{\bf Weaker Dynamical degree comparison conjecture:} Let $f$ be a dynamical correspondence on $X$. Let $N_c, C_c:[0,2n]\rightarrow [0,+\infty )$ be smallest log-concave functions satisfying the following conditions: $N_c(j)\geq \lambda _j(f)$ for all $j=0,1,\ldots, n$, and $C_c(k)\geq \chi _k(f)$ for all $k=0,1,2\ldots ,2n$. Then for all $x\in [0,2d]$, we have $N_c(x)=C_c(x)$.     

Note that Weaker Dynamical degree comparison conjecture is weaker than Dynamical degree comparison conjecture, but they are the same if the sequence $\lambda _j(f)$ is log-concave. The main result of this paper is the following. 

\begin{theorem} Let $X$ be an irreducible smooth algebraic variety of dimension $n$, over $\mathbf{K}$ the algebraic closure of a finite field $\mathbb{F}_q$, and let $H^*(X)$ be a Weil's cohomology theory on $X$. Assume that Standard conjecture D holds for $H^{2n}(X\times X)$. 

1) Then on $X$, Weaker Dynamical degree comparison conjecture holds. 

2) For each $0\leq k\leq 2n$,  there is a positive integer $\rho$ such that if $F$ is an arbitrary polarized endomorphism of $X$, then the largest Jordan block of $F^*:H^k(X)\rightarrow H^k(X)$ has size exactly $\rho$.  

In particular, if there is a polarized endomorphism $F$ of $X$ which is semisimple, then Semisimplicity conjecture holds on $X$. 
 \label{Theorem1}\end{theorem}

Part 2 of Theorem \ref{Theorem1} offers a practically less difficult way to approach the Semisimplicity conjecture than going through all Standard conjectures. Indeed, as mentioned, the Hodge standard conjecture seems more difficult than Standard conjecture D. Moreover, showing semisimplicity of Frobenius morphisms alone may not need proving all the Standard conjectures, for example for the case $X=$ a K3 surface, then semisimplicity for Frobenius morphism was proven by Deligne \cite{Deligne72}.

In particular, we obtain the following consequence. 

\begin{corollary}
If Standard conjecture D holds for $H^{2n}(X\times X)$, then Dynamical degree comparison conjecture holds for all dominant rational maps $f$ on $X$. 
\end{corollary}
\begin{proof}
If $f$ is a dominant rational map on $X$, then its numerical dynamical degrees $\lambda _j(f)$ are a log-concave sequence. Therefore, since Weaker Dynamical degree comparison conjecture is equivalent to Dynamical degree comparison conjecture for dynamical correspondences which have the latter log-concavity property, it follows that Dynamical degree comparison conjecture holds for $f$.  
\end{proof}

Here is a way to construct new examples from known ones for which Theorem \ref{Theorem1} applies. 

\begin{theorem} Let $A\subset X$ be smooth varieties defined over finite fields. We fix a cohomology theory $H^*$. Assume that both $A\times A$, $A\times X$ and  $X\times X$ satisfy Standard conjecture D. Assume moreover that both $A$ and $X$ satisfy the Semisimplicity conjecture. Let $Z=Bl_AX$ be the blowup of $X$ along $A$, and  let $Y$ be a smooth variety so that there is a surjective morphism $\pi :X\rightarrow Y$. Then $Y$ satisfies both the Weaker Dynamical degree comparison conjecture and the Semisimplicity conjecture. Moreover, all rational maps on $Y$ satisfies the Dynamical degree comparison conjecture.  
\label{TheoremNewConstruction}\end{theorem}
An explicit example is the Hilbert scheme $X^{[2]}$ of zero cycles of length $2$ of appropriate varieties $X$ (see Section 2 for more detail).

For proving Theorem \ref{Theorem1}, we combine the ideas in \cite{HT24} to that in \cite{Xie24}. The first main idea in \cite{Xie24} is, over a finite field $\mathbb{F}_q$, to each finite correspondence $f$ and each $k=0,1,\ldots ,n$, consider a composition $Fr_q^s\circ f^t$ where $s$ is an integer (it is important to allow {\bf negative values for $s$}) and $t$ a positive integer. If the dynamical degrees $\lambda _k(f)$ are a log-concave sequence, using that Frobenius morphisms commute with correspondences defined over the same finite field, with a suitable choice of $s$ and $t$ one can make sure that $\lambda _k(Fr_q^s\circ f^t)$ is not too far away compared to $\max _{j=0,\ldots ,n}\lambda _j(Fr_q^s\circ f^t)$. Note that the choice of $\mu _i$ in the proof of Lemma 2.13 in \cite{Xie24} is the same as the choice of $r$ in the proof of Lemma 4.9 in \cite{HT24}. Note that for a finite correspondence $f$, as in the case of endomorphisms, the cohomological dynamical degree $\chi _{2k}(f)$ is the same as the spectral radius of $f^*|_{H^{2k}(X)}$. Therefore, a known result on spectral radii comparison (Theorem 1.2 in \cite{Truong}), which is a generalisation of Entropy comparison conjecture, can be used to show that $\chi _{2k}(f)=\lambda _k(f)$. The inequality $\chi _{2k+1}(f)^2\leq \lambda _k(f)\lambda _{k+1}(f)$ is obtained from this by the trick of using the product correspondence $f\times f:X\times X\rightarrow X\times X$. Thus, Dynamical degree comparison conjecture holds for finite correspondences whose numerical dynamical degrees are a log-concave sequence. In the general case, the numerical dynamical degrees $\lambda _k(f)$ may not be log-concave. The second main idea in \cite{Xie24} is (see Theorem 1.14 therein, and Proposition \ref{Proposition1} below) to establish some interesting relations between certain polytopes of two sequences satisfying certain properties. Theorem 1.14 in \cite{Xie24} can be viewed as a continuous version of Lemma 4.9 in \cite{HT24} (see Proposition \ref{Proposition1} below). 

To prove Theorem \ref{Theorem1}, here we replace the cohomological operator $G_r$ in the approach in \cite{HT24} by composing with iterations of $F^*$. There is a complication arising in case $f$ is not a finite correspondence: in this case the dynamical degree is not the spectral radius but involves norms of dynamical iterates of the correspondence $f^{\Diamond t}$. In case  the dynamical degrees $\lambda _k(f)$ are a log-concave sequence, then again Lemma 4.9 in \cite{HT24} can be used to conclude Theorem \ref{Theorem1}. In the general case, we make use of special properties of the Frobenius norm on matrices and of asymptotic singular values of powers of a matrix (this property was used in \cite{Hu19,Hu24} to solve Dynamical degree comparison conjecture for endomorphisms of Abelian varieties), as well as use essentially the second main idea in \cite{Xie24} mentioned above, but in the following form which is similar to the form of Lemma 4.9 in \cite{HT24}. (More precisely, in Lemma 4.9 in \cite{HT24}, Condition 2 below is required, but instead of Condition 1 below it is required that the sequence $a_j$ is log-concave.)

\begin{proposition} Let $a_0,a_1,\ldots ,a_n$, and $b_0,b_1,\ldots ,b_{2n}$ be non-negative numbers. Assume the following two conditions are satisfied: 

Condition 1: $b_{2j}\geq a_j$ for all $j=0,\ldots, n$. 

Condition 2: For all $r\geq 0$ and for all $k=0,1,\ldots ,2n$, we have: $r^kb_k\leq \max _{j=0,1,\ldots ,n}r^{2j}a_j$.

Let $A_c,B_c:[0,2n]\rightarrow [0,\infty]$ be smallest log-concave functions satisfying the following properties: $A_c(j)\geq a_j$ for all $j=0,1,\ldots ,n$ and $B_c(k)\geq b_k$ for all $k=0,1,\ldots ,2n$.

Then for all $x\in [0,2n]$, we have $A_c(x)=B_c(x)$. 

\label{Proposition1}\end{proposition}

This proposition will be applied to the case $a_j=\lambda _j(f)$ and $b_k=\chi _{k}(f)$, for a dynamical correspondence $f$. In this case, the inequality $a_j\leq b_{2j}$ is always satisfied, since the numerical group $N^j(X)$ is a quotient of the cohomological group $H^{2j}(X)$, see detail in e.g. \cite{Truong20}.

Using the same idea, we also obtain a new proof of the main result in \cite{Xie24}.

\begin{theorem} Let $f:X\vdash X$ be a finite correspondence. Then Weaker Dynamical degree comparison conjecture holds for $f$. If, moreover $f$ is an endomorphism, then Dynamical degree comparison conjecture holds for $f$.  
\label{Theorem2}\end{theorem}

We note that as in the case of an endomorphism, Theorem \ref{Theorem2} holds for any algebraically closed field $\mathbf{K}$, not just the algebraic closure of a finite field. This is done by reducing a general field $\mathbf{K}$ to the algebraic closure of a finite field, via the well known "spreading out" and "specialization arguments". For the case where $f$ is an endomorphism, such an argument is explicitly carried out in \cite{Truong}.

The main results are proven in Section 2. There, we also discuss some heuristic arguments to support Conjecture $G_r$ in \cite{HT24}, by relating it to iterations of polarized endomorphisms. 

{\bf Acknowledgements.} The first author was supported by a grant from Nanjing University (No.~14912209) and a grant from the National Natural Science Foundation of China (No.~12371045). This paper is partly written while the second author is visiting Department of Applied and Computational Mathematics at University of Notre Dame, with a support from the Department of Mathematics, University of Oslo, and he would like to thank both these institutions.

\section{Proof of main results and a discussion on Conjecture $G_r$}
\label{SectionProofs}

In this section we give proofs of Theorems \ref{Theorem1}, \ref{TheoremNewConstruction} and \ref{Theorem2}. We will also discuss a heuristic argument to support Conjecture $G_r$ in \cite{HT24}. Since $X$ is defined over the algebraic closure of a finite field, $X$ itself is defined over a finite field $\mathbf{F}_q$.  We will need the following facts. 

{\bf Fact 1.} As cohomological correspondences, the pullback of Frobenius morphism $Fr_q^*$ commutes with algebraic correspondences, which are defined over $\mathbb{F}_q$, of $X$. 

{\bf Fact 2.} For an endomorphism $F:X\rightarrow X$ and any algebraic correspondence $f:X\vdash X$ which is represented by an effective algebraic cycle and an integer $s$ (maybe negative), the cohomological operator $(F^*)^s\circ f^*:H^*(X)\rightarrow H^*(X)$ is represented by an effective algebraic cycle.

Both facts follow from the following useful lemma. We recall that if $\phi:X\vdash Y$ is a correspondence, then $\phi^\sT:Y\vdash X$ is the correspondence obtained by switching the roles of X and Y. If $\phi$ is algebraic, then also $\phi^\sT$ is algebraic. 

\begin{lemma}[{Lieberman's Lemma, cf.~\cite[Lemma~2.1.3]{MNP13}}]
\label{lemma:Lieberman}
Let $\phi \colon X\vdash Y$, $\psi \colon X'\vdash Y'$, $f \colon X\vdash X'$, and $g \colon Y\vdash Y'$ be correspondences shown as follows:
\[
\xymatrix{
X \ar@{|-}[r]^{f} \ar@{|-}[d]_{\phi} & X' \ar@{|-}[d]_{\psi} \\
Y \ar@{|-}[r]^{g} & Y'. }
\]
Then we have
\begin{align*}
(\phi \times \psi)_*(f) &= \psi\circ f\circ \phi^\sT, \\
(\phi \times \psi)^*(g) &= \psi^\sT\circ g\circ \phi,
\end{align*}
where we think of $\phi \times \psi \in \Z^{2n}(X\times Y\times X'\times Y')_\bQ/\sim$ as a correspondence from $X\times X'$ to $Y\times Y'$ by interchanging the second and the third factors.
In particular, we have the following pullback/pushforward identities on cycle class groups and cohomology groups:
\begin{align*}
((\phi \times \psi)_*(f))_* &= \psi_*\circ f_*\circ \phi^*, \quad ((\phi \times \psi)^*(g))_* = \psi^*\circ g_*\circ \phi_*, \\
((\phi \times \psi)_*(f))^* &= \phi_*\circ f^*\circ \psi^*, \quad ((\phi \times \psi)^*(g))^* = \phi^*\circ g^*\circ \psi_*.
\end{align*}
\end{lemma}

For Fact 1,  we choose $X'=X=Y'=Y$, $\phi=\psi=Fr_q$, the $q$-Frobenius endomorphism and obtain
\[
q^n g = (Fr_q \times Fr_q)^*(g) = Fr_q^\sT\circ g\circ Fr_q.
\]
Note that $Fr_q \circ Fr_q^\sT = q^n \id_X$ (projection formula), note that $ Fr_q^\sT$ simply acts like the pushforward on cohomology groups.
Thus, we get $Fr_q \circ g = g \circ Fr_q$.

For the secon fact, it suffices to take $\phi=F$ and $\psi=\id_X$ in Lieberman's lemma:
\[
(F \times \id_X)_*(f) = f\circ F^\sT \text{ and } (F \times \id_X)^*(g) = g\circ F.
\]

For a linear map $A$ of a finite dimensional vector space over $\mathbf{C}$, we denote by $sp(A)$ the spectral radius of $A$, i.e. the maximum of the set $\{|\lambda |: \lambda$ is an eigenvalue of $A\}$. We recall that $sp(A^t)=sp(A)^t$ for all positive integers $t$.   

\subsection{Proof of Theorem \ref{Theorem2}} We first prove a version of Lemma 4.10 in \cite{HT24}, unconditionally.

\begin{lemma} For any rational number $r>0$ and any $k=0,1,\ldots ,2n$, and for any finite correspondence $f$ we have 
$$r^k|sp (f^*|_{H^k(X)} )|\leq \max _{j=0,\ldots ,n}r^{2j}\lambda _j(f).$$ 

\label{LemmaGeneralized4.10}\end{lemma}
\begin{proof}
Let $s$ an integer (maybe negative) and $t$ a positive integer. Then, since $Fr_q$ commutes with $f$,  by Theorem 1.2 in \cite{Truong}, we have 
\begin{equation}
sp((Fr_q^*)^s \circ (f^{\Diamond t})^*|_{H^k(X)})\leq \max _{j=0,\ldots ,n}\lambda _j((Fr_q^s)\circ (f^{\Diamond t}))\leq C\max _{j=0,\ldots ,n}q^{sj}\lambda   _j(f^{\Diamond t}).
\label{EquationInequalityWithFrobenius}\end{equation}
Also, by Fact 1  and since all eigenvalues of $(Fr_q^*)^s|_{H^k(X)}$ have the same absolute value $|q^{sk/2}|$, we have 
$$sp((Fr_q^*)^s\circ (f^{\Diamond t})^*|_{H^k(X)})=q^{sk/2}sp((f^{\Diamond t})^*|_{H^k(X)}).$$
Since $f$ is a finite correspondence, for all $t$ we have  $(f^{\Diamond t})^*|_{H^k(X)}= (f^*)^t|_{H^k(X)}$, and hence: 
$$sp((Fr_q^*)^s\circ (f^{\Diamond t})^*|_{H^k(X)})= q^{sk/2}sp(f^*|_{H^k(X)})^t.$$

We choose a sequence $s_m$ and $t_m=m$ such that $q^{s_m/(2t_m)}$ converges to $r$. (In particular, choose $s_m=$ the integer part of $2t_m \log _qr$.) Taking the $t_m$-th root on both sides of (\ref{EquationInequalityWithFrobenius}), and then taking limit as $m\rightarrow\infty$, we obtain the desired conclusion.  
\end{proof}

Hence, by using Proposition \ref{Proposition1} (or Lemma 4.9 in \cite{HT24} if the sequence $\lambda _j(f)$ is log-concave) we obtain Theorem \ref{Theorem2}.

\subsection{Proof of Theorem \ref{Theorem1}} We are now ready to give the proof of Theorem  \ref{Theorem1}.

For convenience, for a linear map $A:\mathbf{C}^N\rightarrow \mathbf{C}^N$ we will use the Frobenius norm ({\bf Remark:} Frobenius norm is a general norm on matrices, and it has no relation to Frobenius morphisms), which is defined as follows: 
$$||A||_{Frob}=\sqrt{\sum _{j_1,j_2=1,\ldots ,N}|a_{j_1,j_2}|^2},$$
where $(a_{j_1,j_2})_{j_1,j_2=1,\ldots ,N} $ is the square matrix representing $A$. For such a linear map $A$, we define by $A^{\tau}$ the conjugate transpose of $A$, whose matrix is $(\overline{a}_{j_2,j_1})_{j_1,j_2=1,\ldots ,N}$. If $<.,.>$ is the usual inner product on $\mathbf{C}^N$, then for all $u,v\in \mathbf{C}^N$ we have
$$<Au,v>=<u,A^{\tau}v>.$$
Note that in general $A$ and $A^{\tau}$ do not commute. 

We need the following well known property on Frobenius norm . We include a sketch of proof for convenience. 

\begin{lemma} Let $A$ and $B$ be two square matrices of the same dimension $N$, with coefficients in $\mathbf{C}$. Then, 

$$sp((AA^{\tau})^{-1})^{-1/2}||B||_{Frob}\leq ||AB||_{Frob}\leq sp(AA^{\tau})^{1/2}||B||_{Frob}.$$
In particular, if $sp((AA^T)^{-1})^{-1/2}=sp(AA^{\tau})^{1/2}=\sigma $, then 
$$||AB||_{Frob}=\sigma ||B||_{Frob}.$$

\label{LemmaFrobeniusNormOfProduct}\end{lemma}

\begin{proof}

For each $N\times N$ matrix $C$ with complex coefficients, we have $||C||_{Frob}^2=\Tr (C.C^{\tau})$. Therefore, 
\begin{eqnarray*}
||AB||_{Frob}^2=\Tr [AB(AB)^{\tau}]=\Tr [(AA^{\tau})(BB^{\tau})].
\end{eqnarray*}

Now, applying Inequality II in the paper \cite{PT1970} for the trace of product of two matrices $AA^{\tau}$ and $BB^{\tau}$, which are Hermitian and (semi)positive definite, we obtain what needed. 
\end{proof}

Since $sp(AA^{\tau})^{1/2}$ is a norm on the set of matrices, we obtain the following useful fact,  involved in Lemma \ref{LemmaFrobeniusNormOfProduct}. (The reader can also consult \cite[Theorem 1]{Yamamoto67} and \cite[Lemma 5.4]{Hu19,Hu24}, where more general results concerning other eigenvalues of the matrix $AA^{\tau}$ are given.) 

\begin{lemma} Let $A$ be an $N\times N$ matrix with complex coefficients. Then 
\begin{eqnarray*}
\lim _{m\rightarrow\infty}sp(A^m(A^m)^{\tau})^{1/(2m)}=sp(A). 
\end{eqnarray*}

\label{LemmaSpectralNormAsymptoticSpectralRadius}\end{lemma}

For the remaining of the proof of Theorem \ref{Theorem1}, we assume that Standard Conjecture D holds on $H^{2n}(X\times X)$. We recall the following useful fact (see e.g. \cite{Truong20}): We have a direct sum $H^{2n}(X\times X)=H^{2n}_{alg}(X\times X)\oplus H^{2n}_{tr}(X \times X)$. Here $H^{2n}_{alg}(X\times X)\subset H^{2n}(X\times X)$ is the subspace generated by images of algebraic cycles of $X\times X$  of codimension $2n$, and $H^{2n}_{tr}(X\times X)=\{u\in H^{2n}(X\times X): u.v=0,\empty \forall v\in H^{2n}_{alg}(X\times X)\}$. In particular, the cohomological norm of an algebraic correspondence is equivalent to its numerical norm: $||f^*|_{H^{2n}(X\times X)}||=||f^*|_{N^n(X\times X)}||$.  

Also from now on we fix a polarized endomorphism $F:X\rightarrow X$ (coming together with a positive integer $a>1$ and an ample divisor $L$ such that $f^*(L)\sim aL$).

{\bf A special case:}

If the given polarized endomorphism $F$ is semisimple, then we can choose a convenient inner product of $H^*(X)$ as follows. We choose a basis $\mu _1,\ldots ,\mu _N$ of 
$H^*(X)$ consisting of eigenvectors of $F^*:H^*(X)\rightarrow H^*(X)$. Then we define an inner product as follows: If $u=\sum _ja_j\mu _j$ and $v=\sum _jb_j\mu _j$, then $<u,v>=\sum _ja_j\overline{b_j}$.  

For any correspondence $\Gamma$ of $X$, thought of as a linear map on $H^*(X)$, we then define $\Gamma ^{\tau}$ by using the given inner product $<.,.>$. Note that then the matrix of $F^*$ is the diagonal matrix consisting of eigenvalues of $F^*$, and hence $(F^*)^{\tau }$ is the diagonal matrix consisting of conjugates of eigenvalues of $F^*$. 

Now for $k=0,1,\ldots ,2n$, we consider $A=F^*|_{H^k(X)}$. Then from the above paragraph and the fact (see Theorem \ref{Theorem2}) that all eigenvalues of $A$ have eigenvalues of absolute value $a^{k/2}$, we obtain then 
$$sp(AA^{\tau})^{1/2}=sp((AA^{\tau})^{-1})^{-1/2}=a^{k/2}.$$
Therefore, if $f:X\rightarrow X$ is an algebraic correspondence, we obtain from Lemma \ref{LemmaFrobeniusNormOfProduct} the following equality
\begin{equation}
||F^*\circ f^*|_{H^{k}(X)}||_{Frob}=a^{k/2}||f^*|_{H^{k}(X)}||_{Frob}.
\label{EquationEqualityFrobeniusNorm}\end{equation}

\underline{\bf Proof of part 1:}

For each integer $s$ (both positive and negative are allowed), and a positive integer $t$, we consider the cohomological operator $(F^*)^s\circ (f^{\Diamond t})^*$, where $f^{\Diamond t}$ is the t-th dynamical iteration of $f$.  
By definition, $f^{\Diamond t}$ is an effective cycle. Hence, by Fact 2, $(Fr^*)^s\circ (f^{\Diamond t})^*$ is represented by an effective cycle. 

Now, we fix any $k=0,1,\ldots, 2\dim (X)-1, 2\dim (X)$. We will prove the following:  

For every real number $r\geq 0$, we have 
\begin{equation}
r^k\chi _k(f)\leq \max _{j=0,\ldots ,\dim (X)}r^{2j}\lambda _j(f).
\label{Equation1}\end{equation} 
By continuity, it is enough to consider the case $r\not=1$ and $r>0$. 

To this end, we first prove the next result is a strengthen of Lemma 4.10 in \cite{HT24}, here Conjecture $G_r$ is not needed. We define $A=F^*|_{H^k}(X)$. 

\begin{lemma} 
Assume Standard conjecture D holds for $H^{2n}(X\times X)$. There is a constant $C>0$, independent of $f$ (a dynamical correspondence), $s$ (an integer) and $t$ (a positive integer), so that 
$$sp([A^s(A^s)^{\tau}]^{-1})^{-1/2}||(f^{\Diamond t})^*|_{H^k(X)}||_{Frob}\leq C\max _{j=0,\ldots ,\dim (X)}a^{sj}||(f^{\Diamond t})^*|_{N^j(X)}||.$$
\label{LemmaClaim1}\end{lemma}
\begin{proof}
By Theorem \ref{Theorem2}, all eigenvalues of $(F^*)^s$ on $H^k(X)$ have absolute value $a^{sk/2}$. Therefore, by Lemma \ref{LemmaFrobeniusNormOfProduct} (see also (\ref{EquationEqualityFrobeniusNorm}, with $F^*$ replaced by $(F^*)^s$)) we have 
$$sp([A^s(A^s)^{\tau}]^{-1})^{-1/2}||(f^{\Diamond t})^*|_{H^k(X)}||_{Frob}\leq ||(F^*)^s\circ (f^{\Diamond t})^*|_{H^k(X)}||_{Frob}.$$ 

Now, if we define $A_{s,t}=(F^*)^s\circ (f^{\Diamond t})^*|_{H^k(X)}$, then $A_{r,s}$ is algebraic and effective. Note that $A_{r,s}^{\tau}$ may not be algebraic.  We denote by $B_{r,s}$ the image of $A_{r,s}^{\tau}$ in the summand $H^{2n}_{alg}(X\times X)$ of $H^{2n}(X\times X)$.

We have, $$||(F^*)^s\circ (f^{\Diamond t})^*|_{H^k(X)}||_{Frob}^2=\Tr (A_{s,t}A_{s,t}^{\tau}),$$
which by Lefschetz trace formula (here $\Delta _{2n-k}$ is the summand in $H^{2n-k}(X)\otimes H^{k}(X)$ of the diagonal $\Delta$, which is an algebraic cycle by Standard conjecture C, the latter being a consequence of Deligne's proof of Weil's Riemann hypothesis, see \cite{KM74})
\begin{eqnarray*}
|\Tr (A_{s,t}A_{s,t}^{\tau})|&=&|A_{s,t}A_{s,t}^{\tau}.\Delta _{2n-k}| \leq C_1 ||\Delta _{2n-k}|_{H^{2n}(X\times X)}||\times ||A_{s,t}|_{H^{2n)(X\times X)}}||\times ||A_{s,t}^{\tau}|_{H^{2n)(X\times X)}}||\\
&\leq& C_2  ||A_{s,t}|_{H^{2n}(X\times X)}||^2\leq C_3||A_{s,t}|_{N^n(X\times X)}||^2. 
\end{eqnarray*}
Here, the first inequality follows since $A_{s,t}A_{s,t}^{\tau}.\Delta _{2n-k}$ is a trilinear functional on a finite dimensional vector space. The second inequality follows from the fact that the (Frobenius) norm of a linear map and its conjugate transpose are the same, and the norm of $\Delta _{2n-k}$ is a positive constant. The third inequality follows from the above mentioned fact that if we assume Standard conjecture D, then the cohomological norm of an algebraic cycle is equivalent to its numerical norm.  

Since $F$ is polarized, for any $u\in N^j(X)$ with $cL^j-u$ is effective, we have (see \cite{HT24}) $ca^jL^j=c(F^*L^j)\geq cF^*(u)$. Therefore, combining with the fact that all eigenvalues of $F^*$ on $H^{2j}(X)$ have absolute value $a^j$, we obtain that $||F^*|_{N^j(X)}||\sim a^j$. Similarly, for any integer $s$, we have $||(F^*)^s|_{N^j(X)}||\sim a^{sj}$.

Now, since $A_{s,t}$ is represented by an effective algebraic cycle, we have for a fixed ample divisor $H$ on $X$:  
$$||A_{s,t}|_{N^n(X\times X)}||\sim \sum _{j=0}^nA_{s,t}H^j.H^{n-j}=\sum _{j=0}^n(F^*)^s((f^{\Diamond t})^*H^j).H^{n-j}\leq C_3\max _{j=0,\ldots ,n} a^{sj}||(f^{\Diamond t})^*|_{N^j(X)}||.$$

Combining all the above estimates, we obtain the conclusion. 
\end{proof}

Now, we finish the proof of (\ref{Equation1}). We choose a sequence $s_m,t_m=m$ where $s_m$ is an integer (maybe positive or negative) and $t_m=m$ converges to $\infty$, so that $a_1^{s_m/(2t_m)}\rightarrow r$. A way to choose is $t_m=m$ and $s_m=[2m \log _{a_1}r]$ the integer part of $2m \log _{a_1}r$. Since here we assume that $r\not =1$, we have $\lim _{m\rightarrow\infty}|s_m|=\infty$. Applying Lemma  \ref{LemmaClaim1} to for each of the pair $s_m,t_m=m$, we obtain: 
\begin{eqnarray*}
sp([A^{s_m}(A^{s_m})^{\tau}]^{-1})^{-1/2}||(f^{\Diamond m})^*|_{H^k(X)}||_{Frob}\leq C\max _{j=0,\ldots ,\dim (X)}a_1^{s_mj}||(f^{\Diamond m})^*|_{N^j(X)}||. 
\end{eqnarray*}
Taking the  $m$-th root of both sides of the previous inequality, and then limsup as $m\rightarrow \infty$, we obtain: 

$$\limsup _{m\rightarrow\infty}[sp([A^{s_m}(A^{s_m})^{\tau}]^{-1})^{-1/2}||(f^{\Diamond m})^*|_{H^k(X)}||_{Frob}]^{1/m}\leq \max _{j=0,\ldots ,\dim (X)}r^{2j}\lambda _j(f).$$
By Lemma \ref{LemmaSpectralNormAsymptoticSpectralRadius}, we have 
\begin{eqnarray*}
\lim _{m\rightarrow\infty} [sp([A^{s_m}(A^{s_m})^{\tau}]^{-1})^{-1/2}]^{1/m}=[\lim _{m\rightarrow\infty} [sp([A^{s_m}(A^{s_m})^{\tau}]^{-1})^{-1/2}]^{1/(2s_m)}]^{s_m/m}=r^{k}. 
\end{eqnarray*}
Here, in the last equality, we again use the fact from Theorem \ref{Theorem2} that all eigenvalues of $F^*|_{H^k(X)}$ have absolute value $a_1^{k/2}$.

Hence, by using Proposition \ref{Proposition1} (or Lemma 4.9 in \cite{HT24} if the sequence $\lambda _j(f)$ is log-concave), we have that Weaker Dynamical degree comparison conjecture holds.

\underline{\bf Proof of part 2:} 

Recall that $F$ is the given polarized endomorphism, with $F^*(L)=aL$ for a positive integer $a>1$ and an ampled divisor $L$. We fix $0\leq k\leq 2n$, and let $b$ be the largest size of the Jordan blocks of $F^*:H^{k}(X)\rightarrow H^{k}(X)$. Then for integers $s$ we have $||(F^*)^s|_{H^k(X)}||\sim |s|^{b-1}a^{s}$ (this can be seen by applying the function $z\mapsto z^s$ - which is analytic where $s\not= 0$- to the Jordan blocks).

Let now $f:X\rightarrow X$ be a polarized endomorphism (coming together with an ample divisor  $L_f$ and a positive integer $a_1>1$ such that  $f^*(L)\sim a_1L_f$). Then it follows that there is a constant $C>0$ so that for all positive integer $t$ and all $j=0,1,\ldots ,n$ we have $||(f^t)^*|_{N^j(X)}||\leq Ca_1^{jt}$. Moreover, since $f$ is a polarized endomorphism, by Theorem \ref{Theorem2}, all eigenvalues of $f^*|_{H^k(X)}$ have absolute value $a_1^{k/2}$. Therefore, if $b_1$ is the largest size of the Jordan blocks of $f^*:H^k(X)\rightarrow H^k(X)$, then for positive integers $t$ we have $||(f^*)^t|_{H^k(X)}||\sim t^{b_1-1}a_1^{t}$.

We want to show that  $b_1=b$. By symmetry, it suffices to show that $b_1\leq b$. We assume otherwise that $b_1>b$, and will arrive at a contradiction. 

Lemma \ref{LemmaClaim1} then yields that for all integers $s,t$: 
$$\frac{a^{sk/2}}{|s|^{b-1}}||(f^t)^*|_{H^k(X)}||_{Frob}\leq C \max _{j=0,1,\ldots ,n}a^{sj}a_1^{tj}.$$

From the above we obtain the inequality: 
$$a^{sk/2}a_1^{tk/2}\frac{t^{b_1-1}}{|s|^{b-1}}\leq C\max _{j=0,1,\ldots ,n}a^{sj}a_1^{tj},$$
for all $s$ integer and $t$ positive integer. Put $r=a^{s}a_1^{t}$, we then get 
\begin{equation}
r^{k/2}\frac{t^{b_1-1}}{|s|^{b-1}}\leq C\max _{j=0,1,\ldots ,n}r^j.
\label{EquationContradiction}\end{equation}
We consider several cases, depending on the values of $a$ and $a_1$. We can write $r=a^sa_1^t=e^{\log a_1 (\frac{\log a}{\log a_1}s +t )}$. 

Case 1: $\frac{\log a}{\log a_1}\in \mathbf{Q}$. In this case there is a sequence $s_m$ and $t_m\rightarrow \infty$ so that $a^{s_m}a_1^{t_m}=1$ for all $m$. Note that $|t_m/s_m|$ is bigger than a constant $c>0$ uniformly in $m$. Then we obtain a contradiction $t_m^{b_1-1}/|s_m|^{b-1}\leq C_1$ for all $m$ if $b_1>b$. 

Case 2: $\frac{\log a}{\log a_1}\notin \mathbf{Q}$. In this case, by Kronecker's approximation theorem \cite{Kronecker1884}, there are integers $s_m$ and $t_m>0$, so that 
$$|\frac{\log a}{\log a_1}s_m +t_m |<\frac{1}{m}.$$
It follows easily that $t_m $ converges to $+\infty$, and moreover $|t_m/s_m|$ is bigger than a constant $c>0$ uniformly in $m$. Again from (\ref{EquationContradiction}), we obtain a contradiction
$$e^{-\frac{k|\log a_1|}{2m}}\frac{t_m^{b_1-1}}{|s_m|^{b-1}}\leq C\max _{j=0,1,\ldots ,n}e^{\frac{j|\log a_1|}{m}},$$
for all $m$.

\subsection{Proof of Theorem \ref{TheoremNewConstruction} and some explicit examples} We first start with a general result on the relation between blowups and Standard conjecture D and Semisimplicity conjecture. 

\begin{proposition} Let $A\subset X$ be smooth varieties defined over finite fields, and $Z=Bl_A(X)$ the blowup of $X$ along $A$. 

1) If both $A$ and $X$ satisfies Standard conjecture D, then $Z$ also satisfies Standard conjecture $D$. 

2) If the Frobenius morphisms $Fr_A$ and $Fr_X$ are both semisimple, then the morephism $Fr_Z$ is also semisimple.

3) If both $A\times A$, $A\times X$ and $X\times X$ satisfy Standard conjecture $D$, then $Z\times Z$ also satisfies Standard conjecture D.

4) Let $\pi :Z\rightarrow Y$ be a surjective morphism, where $Y$ is also smooth. If $Z$ satisfies Standard conjecture D then $Y$ also satisfies Standard conjecture D. If $Fr_Z$ is semisimple, then $Fr_Y$ is also semisimple. 

\label{PropositionBlowupAndConjectures}\end{proposition}

\begin{proof}

 These properties are probably well-known to experts. For example, in \cite{Arapura2006}, it is shown that over $\mathbf{C}$, the Lefschetz standard conjecture holds on $Z$ if it holds for both $A$ and $X$. Here we provide a sketch of proof for the reader's convenience. Let $\pi :Z\rightarrow X$ be the blowup map, and $E=\pi ^{-1}(A)$ the exceptional divisor.

Let $r$ be the codimension of $A$ in $X$. We recall that the cohomology group of $Z$ is related to that of $X$ and $A$ by the following relation: 
\begin{eqnarray*}
H^k(Z)=\pi ^*H^k(X)\oplus V_k, 
\end{eqnarray*}
where 
\begin{eqnarray*}
V_k=\oplus _{0\leq j\leq r-1}(E^j.\pi ^*H^{k-2j}(A)). 
\end{eqnarray*}

In particular, assume that $\alpha \in H^k(Z)$ is such that $\pi _*(\alpha )=0$. Then $\alpha $ belongs to $V_k$. Moreover, $\alpha$ is $0$ iff for all $\beta \in V_{2n-k}$ then $\alpha .\beta =0$. (Indeed, by Poincar\'e duality, $\alpha =0$ iff $\alpha .\beta =0$ for all $\beta \in H^{2n-k}(X)$. Since $\pi _*(\alpha )=0$, if $\beta \in \pi ^*H^{2n-k}(X)$ then $\alpha .\beta =0$, by projection formula. Here we use that if $0\leq j\leq r-1$, then since $E^j$ has dimension at least $n-r+1$ and $\pi _*(E^j)$ is an algebraic cycle with support in $A$, which has dimension $n-r$, we must have $\pi _*(E^j)=0$.)

1) Let $B$ be an algebraic cycle on $X$. Then we can write (in the Chow ring) $B=\pi ^*\pi _*(B)+B_1$, where $B_1$ is also an algebraic cycle so that $\pi _*(B_1)=0$.

If $B$ is numerically $0$ in $Z$, then $\pi _*(B)$ is numerically $0$ in $X$ and $B_1$ is numerically $0$ in $Z$. 

Since $X$ satisfies Standard conjecture D, it follows that $\pi _*(B)$ is homologically $0$. Hence, also $\pi ^*\pi _*(B)$ is homologically $0$ in $Z$. Therefore, to show that $B$ is homologically $0$, we only need to show that $B_1$ is homologically $0$. Since $B_1$ is numerically $0$ and $\pi _*(B_1)=0$, we can see (from the discussion before part 1) that this is reduced to the assumption that $A$ satisfies Standard conjecture D. 

2) By Poincar\'e duality, we note that $Fr_Z$ is semisimple (meaning $Fr_Z^*:H^*(X)\rightarrow H^*(X)$ is semisimple) iff $(Fr_Z)_*: H^k(X)\rightarrow H^k(X)$ is semisimple. By the decomposition mentioned before part 1, we have $(Fr_Z)_*:\pi^*H^k(X)\rightarrow \pi^*H^k(X)$ and $(Fr_X)_*:V_k\rightarrow V_k$. Now, $(Fr_Z)_*:\pi ^*H^*(X)\rightarrow \pi^*H^*(X)$ is semisimple is more or less the same as that $Fr_X$ is semisimple, and   $(Fr_X)_*:V_k\rightarrow V_k$ is semisimple is more or less the same as that $Fr_A$ is semisimple. The latter are assumed, hence we are done. 

3) This follows from part 1 by observing that: 
\begin{eqnarray*}
Bl_AX\times Bl_AX&=&Bl_{A\times Bl_AX}(X\times Bl_AX),\\
A\times Bl_AX&=&Bl_{A\times A}(A\times X),\\
X\times Bl_AX&=&Bl_{X\times A}(X\times X).
\end{eqnarray*}
\end{proof}

Now we are ready to prove Theorem \ref{TheoremNewConstruction}. 
\begin{proof}[Proof of Theorem \ref{TheoremNewConstruction}] Let $Z=Bl_AX$. Since both $A\times A$, $A\times X$ and $X\times X$ satisfy Standard conjecture D, by Proposition \ref{PropositionBlowupAndConjectures}, $Z\times Z$ also satisfies Standard conjecture D. 

Since both $A$ and $X$ satisfy semisimplicity conjecture, both $Fr_A $ and $Fr_X$ are semisimple. Then, by  Proposition \ref{PropositionBlowupAndConjectures} again, $Fr_Z$ is semisimple.

Since $Z\times Z$ satisfy Standard conjecture D and $Fr_Z$ is semisimple, by Theorem \ref{Theorem1} we have that the Semisimplicity conjecture holds on $Z$. Also, the Weaker Dynamical degree comparison conjecture holds for $Z$. 

Now, since $\pi :Z\rightarrow Y$ is a surjective morphism, Proposition \ref{PropositionBlowupAndConjectures} implies that $Y\times Y$ satisfies Standard conjecture D and $Fr_Y$ is semisimple. Then by Theorem \ref{Theorem1}, we have that $Y$ satisfies the Semisimplicity conjecture, and also the Weaker Dynamical degree comparison conjecture.  

4) This comes easily from the fact that the pullbacks by $\pi ^*$ of both numerical and cohomological groups are injective, and the Frobenius morphism commutes with $\pi$.  
\end{proof}

For an example, let $X$ be such that $Fr_X$ is semisimple and $X^4$ satisfies the Standard conjecture D. Let $X^{[2]}$ be the Hilbert scheme of zero cycles of length 2 on $X$. Then $X^{[2]}$ is smooth, and more  over $X^{[2]}=Bl_{\Delta }(X\times X)/S_2$, where $\Delta \subset X\times X$ is the diagonal, and $S_2$ is the symmetric group of 2 elements. Since $\Delta $ is isomorphic to $X$, assumptions in Theorem \ref{TheoremNewConstruction} are satisfied (note that $Fr_{X\times X}=Fr_X\times Fr_X$ is semisimple if $Fr_X$ is semisimple). We conclude that $X^{[2]}$ satisfies the Semisimplicity conjecture, and also the Weaker Dynamical degree comparison conjecture. 

For an explicit $X$ which satisfies the conditions in the previous paragraph, we can choose $X$ an Abelian variety \cite{Clozel99}, certain self product of $K3$ surfaces \cite{IIK22}, or any multi-projective space $\mathbf{P}^{n_1}\times \ldots \times \mathbf{P}^{n_k}$ (simple use of Kunneth's formula).   

\subsection{A discussion on Conjecture $G_r$} We first recall Conjecture $G_r$ from \cite{HT24}. 

{\bf Conjecture $G_r$.} For any $r\in \bQ_{>0}$, the cohomological correspondence $G_{r}$ of $X$ is algebraic and represented by a rational algebraic $n$-cycle - still denoted  $G_r$ - on $X\times X$.
Moreover, for any dynamical correspondence $f$ of $X$, there is a constant $C>0$, independent of $r$ and $f$, such that
\begin{equation}
\label{eq:Gr}
\|G_r\circ f\| \le C \deg(G_r\circ f),
\end{equation}
where $\norm{\cdot}$ denotes a norm on $\N^n(X\times X)_\bR$ and $\deg(\cdot)$ is the degree function on $X\times X$ with respect to a fixed ample divisor $H_{X\times X}$.

Since $\mathbf{K}$ is the algebraic closure of a finite field, the first part of the conjecture (i.e. that $G_r$ is algebraic) follows from Weil's Riemann hypothesis. Therefore, it remains to prove the second part, i.e. the inequality $\|G_r\circ f\| \le C \deg(G_r\circ f)$. 

Note that if $G_r$ is represented by an {\bf effective} algebraic cycle, then the concerned inequality is easily proven. When $X$ is an Abelian variety, then $G_r$ can  be related to the usual multiplication maps on $X$, whose graphs are in particular effective algebraic cycles, and hence Conjecture $G_r$ is proven, see \cite{HT24}. 

In general, we do not know if $G_r$ can be represented by an effective algebraic cycle. However, the Frobenius morphisms behave similar to $G_r$ in many aspects. First, by the ideas in \cite{Xie24}, an iterate of the Frobenius morphism $(Fr_q^*)^s$ acts as multiplication by $r^{k/2}$ on $sp(f^*|_{H^k(X)})$ for $r=q^{1/2}$. Second, by the ideas in the proof of Theorem \ref{Theorem1}, if we assume that there is a polarized endomorphism $F$   which is semisimple, then  $(F^*)^s$ also acts as multiplication by $r^{k/2}$ on $||f^*|_{H^k(X)}||_{Frob}$, for an appropriate $r>0$.


\bibliographystyle{amsalpha}

\bibliography{mybib}

\end{document}